\documentclass[12pt]{article}%
\usepackage{amsfonts}
\usepackage{amsmath}
\usepackage{amssymb}
\usepackage{graphicx}
\usepackage{subfigure}
\usepackage{citesort}

\newtheorem{theorem}{Theorem}

\begin{document}

\author{Leonid Kunyansky }
\title{Reconstruction of a function from its spherical (circular) means with the
centers lying on the surface of certain polygons and polyhedra}
\maketitle

\begin{abstract}
We present explicit filtration/backprojection-type formulae for the
inversion of the spherical (circular) mean transform with the centers lying on
the boundary of some polyhedra (or polygons, in 2D). The formulae are derived
using the double layer potentials for the wave equation, for the domains with
certain symmetries.  The formulae are valid for a rectangle and certain
triangles in 2D, and for a cuboid, certain right prisms and a certain pyramid
in 3D. All the present inversion formulae yield exact reconstruction within
the domain surrounded by the acquisition surface even in the presence of
exterior sources.
\end{abstract}


\section*{Introduction}

The spherical mean Radon transform and related inversion formulae are
important, in particular, for solving problems of thermo- and photo-acoustic
tomography (TAT/PAT) \cite{kruger1,kruger,Oraev94}. In these modalities a
region of interest (ROI)\ within a human body is subjected to a very short
electromagnetic (EM) pulse which causes thermoelastic expansion of the tissues
and generates an outgoing acoustic wave. The resulting acoustic pressure is
measured by the detectors placed on a surface surrounding ROI. The initial
pressure distribution $f(x)$ is not uniform: tumors absorb much more
EM\ energy than healthy tissue and thus generate a much stronger signal.
Hence, by recovering $f(x)$ one obtains valuable medical information.

Under certain simplifying assumptions the problem of finding the initial
pressure can be re-stated as that of inverting the spherical mean Radon
transform in 3D over spheres with centers lying at the measuring surface. A
similar problem in 2D arises when instead of usual point-like detectors one
measures the acoustic pressure with linear detectors
\cite{Burgh,haltmeierinterf,Paltaufnew1}. A detailed discussion of the
solvability, range conditions, and inversion techniques for the spherical mean
transform in 2D and 3D can be found in
\cite{AK,AKK,AKQ,AQ,FPR,FR,FR2,kuchre,KuKu,Kukuspringer,kunya,PatchSch,WangCRC,Wangbook,MXWreview}
and references therein.

In the present paper we derive explicit inversion formulae for the spherical
(or, in 2D, circular) mean transform for certain measurement surfaces with
special symmetries. In general, such explicit formulae play an important role
in inverse problems: they provide an important theoretical insight and, in
some cases, they serve as a starting point for the development of efficient
reconstruction algorithms. For example, a very popular
filtration/backprojection (FBP) algorithm (see, for example, \cite{Natt}) is
obtained as a result of discretization of one of the explicit inversion
formulae for the 2D Radon transform.

For the spherical mean Radon transform the inversion formulae of the
filtration/back\-projection type are known only for special types of the
acquisition surfaces. In particular, several different formulae for the
spherical measuring surface have been found in \cite{XuWang0,XuWang} (in 3D),
in \cite{FPR} (in odd dimensions), in \cite{FHR} (in even dimensions), and in
\cite{kunya} (in arbitrary dimensions). A general family containing as
particular cases all the above-mentioned formulae for a sphere, was derived in
\cite{nguyen}. The so-called "universal" backprojection formula~\cite{XuWang}
(already mentioned above) holds not only for a spherical measuring surface,
but also for a the surface of an infinite cylinder and for a plane
\cite{XuWang1,XuWang} (both in 3D). Other formulae for the planar acquisition
are also known \cite{norton,Palamodov,Fawcett,denis,Burgh}. The existence of
the FBP-type formulae for these three types of surfaces (and the absence of
such formulae for other smooth surfaces) can be explained by certain
symmetries common to a sphere, a plane, and a cylinder, as discussed in
\cite{Kukuspringer}.

For other measuring surfaces there exist certain numerical techniques for
solving the problems TAT/PAT (see
\cite{AK,Kukuspringer,PatchSch,WangCRC,Wangbook,MXWreview} and references
therein); but the inversion formulae for such surfaces are not known. (The
notable exception is the explicit general solution obtained in \cite{AK} for
any closed surface; however, this solution is given in terms of functions of
the wave operator rather than in a form of a closed form integro-differential
expression). In the present paper we derive explicit FBP-type inversion
formulae for the spherical mean transform with the centers lying on the
surfaces of certain polyhedra (in 3D) or on the boundaries of certain polygons
(for circular mean transform in 2D).

The paper is organized as follows. In Section~\ref{S:formula} we provide a
formal description of the problem and supply some introductory information.
Next, in Section~\ref{S:wave} we explicitly solve the time reversal problem
for the wave equation in a cube or a square (in 2D) by means of the double
layer wave potentials. Such a solution is possible due to certain symmetries
of these regions. Next, in Section~\ref{S:means}, from the explicit solution
of the wave equation we derive FBP-type formulae for the inversion of the
spherical (circular) mean transform with centers on the boundary of a cube (a
square). The present inversion formulae have an interesting property of being
insensitive to acoustic sources lying outside of the region surrounded by the
acquisition surface, as discussed in Section~\ref{S:exterior}. Numerical
examples of image reconstructions obtained using our formulae are provided in
Section~\ref{S:numerics}. Finally, in Section~\ref{S:others} we apply the same
methods to some other domains possessing the necessary symmetries. In 2D these
include a rectangle and three types of triangles: an equilateral and a right
isosceles triangles, and a triangle with the angles $\frac{\pi}{2},\frac{\pi
}{3}$ and $\frac{\pi}{6}$. In 3D our formulae extend to cuboids, right prisms
whose base is one of the three triangles mentioned above, and a pyramid whose
side faces are equal right isosceles triangles.

\section{Formulation of the problem\label{S:formula}}

Throughout the paper we assume that the source of the acoustic pressure is
supported within a bounded region $\Omega$, the detectors are placed on the
boundary $\partial\Omega,$ and the speed of sound is constant and equal to $1$
in the whole space $\mathbb{R}^{n}.$ The time-dependent acoustic pressure
$p(x,t)$ satisfies the wave equation
\begin{equation}%
\begin{cases}
p_{tt}=\Delta_{x}p,\quad t\geq0,\quad x\in\mathbb{R}^{n}\\
p(x,0)=f(x),\quad p_{t}(x,0)=0,
\end{cases}
\label{E:wave_eq}%
\end{equation}
where initial pressure $f(x)$ is supported within $\Omega.$ The detectors
measure the values $P(y,t)$ of the pressure on the boundary:%
\begin{equation}
{P(y,t)=}\left.  {p(y,t)}\right\vert _{{y\in\partial\Omega}}{,\quad}%
t\in\lbrack0,\infty){.} \label{E:wavedata}%
\end{equation}
Our goal is to reconstruct $f(x)$ from measurements ${P(y,t).}$

When the measurements are done using conventional point-like acoustic
detectors, the above problem has to be solved in 3D setting. However, when the
so-called line detectors are utilized \cite{Burgh,haltmeierinterf,Paltaufnew1}%
, the 3D problem reduces to a set of similar 2D inverse problems for the wave
equation, so from the practical point of view both 2D and 3D cases are of interest.

In the 3D case the measurements need only be done during the time interval
$[0,\mathrm{diam}\Omega],$ since the signal vanishes after $t=\mathrm{diam}%
\Omega$ due to the Huygens principle \cite{Vladimirov}. In the 2D case, the
values of $P(y,t)$ in the whole interval $t\in\lbrack0,\infty)$ are needed for
the theoretically exact reconstruction, and we will assume that these values
are known. (However, if $P(y,t)$ is known for all
$t\in\lbrack0,\mathrm{diam}\Omega]$, the remaining values for the time
interval $[\mathrm{diam}\Omega,\infty)$ can be explicitly reconstructed, see \cite{PaltCRC}
and the discussion at the end of the next section).

\subsection{Spherical means}

The problem of finding $f(x)$ can also be re-formulated in terms of inverting
the spherical mean transform. Indeed, in 3D the solution to the initial value
problem (\ref{E:wave_eq}) can be written in the form of the Kirchhoff formula%
\begin{equation}
p(x,t)=\frac{\partial}{\partial t}tM(x,t),\label{E:3daverages}%
\end{equation}
where $M(x,r)$ are the spherical means of $f(x)$%
\[
M(x,r)=\frac{1}{4\pi}\int\limits_{S^{2}}f(x+r\omega)d\omega,
\]
and where $S^{2}$ is the unit sphere in 3D. Thus, using (\ref{E:wavedata}) one
can relate the spherical means $M(y,t)$ to ${P(y,t)}$:%
\[
{P(y,t)=}\frac{\partial}{\partial t}tM(y,t),
\]
and, taking into account that $M(y,0)=0,$ solve for $M(y,t)$%
\[
M(y,t)=\frac{1}{t}\int_{0}^{t}P(y,t^{\prime})dt^{\prime}.
\]
Now the problem is reduced to finding function $f(x)$ from its spherical means
$M(y,t)$ with centers on $\partial\Omega.$

Similarly, in 2D, one can introduce the circular means $m(x,r):$%
\[
m(x,r)=\frac{1}{2\pi}\int\limits_{S^{1}}f(x+r\omega)d\omega.
\]
Now solution of (\ref{E:wave_eq}) can be expressed through $m(x,r)$ by the
formula%
\[
p(x,t)=\frac{\partial}{\partial t}\int\limits_{0}^{t}\frac{m(x,r)rdr}%
{\sqrt{t^{2}-r^{2}}},
\]
and, in particular,%
\begin{equation}
P(y,t)=\frac{\partial}{\partial t}\int\limits_{0}^{t}\frac{m(y,r)rdr}%
{\sqrt{t^{2}-r^{2}}},y\in\partial\Omega. \label{E:Abel}%
\end{equation}
The integral in the right hand side of the latter equation is the well known
Abel transform whose inversion formula is also well known. Thus, circular
means $m(y,r)$ can be found from $P(y,t)$ by the following formula%
\[
m(y,r)=\frac{2}{\pi}\int\limits_{0}^{r}\frac{P(y,t)dt}{\sqrt{r^{2}-t^{2}}},
\]
which again leads to the problem of reconstructing $f(x)$ from its cylindrical
means $m(y,r).$ (Since $m(y,r)$ vanish for $r>\mathrm{diam}\Omega,$ the
knowledge of $P(y,t)$ in the time interval $t\in\lbrack0,\mathrm{diam}\Omega]$
is sufficient to reconstruct $m(y,r)$ using the latter formula; moreover,
values of $P(y,t)$ for $t>\mathrm{diam}\Omega$ can now be recovered from
$m(y,r)$ by equation (\ref{E:Abel}).)

\subsection{Reconstruction strategies}

It is well known that the initial pressure $f(x)$ can be found by the time
reversal, i.e. by solving the wave equation back in time. In particular, in
the 3D case, since solution to (\ref{E:wave_eq}) vanishes within $\Omega$
after $t=\mathrm{diam}\Omega$ due to the Huygens principle, one can solve the
following initial-boundary-value problem backwards in time :%
\begin{equation}%
\begin{cases}
u_{tt}=\Delta_{x}u,\quad t\in\lbrack0,T],\quad x\in\Omega\subset\mathbb{R}%
^{3},\\
u(x,T)=0,\quad u_{t}(x,T)=0,\\
u(y,t)=P(y,t),\quad y\in\partial\Omega,\\
T=\mathrm{diam}\Omega.
\end{cases}
\label{E:3dbackwards}%
\end{equation}

Function $f(x)$ we seek to recover equals $u$ at the moment $0$:%
\begin{equation}
f(x)=u(y,0). \label{E:terminal}%
\end{equation}

In 2D the solution to the wave equation does not vanish in a finite time, and
therefore the initial condition should be replaced by the condition of the
decrease at $t\rightarrow\infty:$
\begin{equation}%
\begin{cases}
u_{tt}=\Delta_{x}u,\quad t\in\lbrack0,\infty),\quad x\in\Omega\subset
\mathbb{R}^{2},\\
\lim\limits_{T\rightarrow\infty}u(x,T)=0,\quad\lim\limits_{T\rightarrow\infty
}u_{t}(x,T)=0,\\
u(y,t)=P(y,t),\quad y\in\partial\Omega.
\end{cases}
\label{E:back2d}%
\end{equation}
As in 3D, $f(x)$ is found using formula (\ref{E:terminal}).

The above initial-boundary-value problems can be solved numerically using, for
example, finite difference methods (e.g. \cite{AmbPatch,Burgh1}). For certain
acquisition surfaces methods based on eigenfunction expansion can also be very
effective \cite{kunyaser}. Our goal, however, is to obtain an explicit
filtration/backprojection type formula for the solution. If the existing
detectors could measure (in addition to the pressure) the normal derivative of
the pressure, such an explicit analytic solution could have been represented,
using the Green's formula, in the form of the double- and single- layer
potentials. Unfortunately, the normal derivative is not known. Still, explicit
representations of the solution of the wave equation by means of the double
layer potential are possible for surfaces with certain symmetries, such as,
for example, a sphere, a cylinder, and a plane. In fact, the so-called
"universal backprojection formula" \cite{XuWang} for these surfaces can be
interpreted as evaluation of a certain double layer potential.

In the next section we find explicit inversion formulae for the case when
$\Omega$ is a cube (in 3D) or a square (in 2D). First, we obtain explicit
solutions for the problems (\ref{E:3dbackwards}) and (\ref{E:back2d}) in the
form of the double layer potentials. These formulae allow us to reconstruct
$f(x)$ from the measured values of the pressure $P(y,t);$ in
Section~\ref{S:means} they also will be used to derive explicit
filtration/backprojection time formulae for recovering $f(x)$ from its
spherical (circular) means in these domains.

\section{Explicit solution of the wave equation in a cube or a
square\label{S:wave}}

\subsection{Double layer potentials}

In 3D, the retarded free space Green's function $G^{3D+}(x,t)$ is given by the
following formula~\cite{Vladimirov}%
\[
G^{3D+}(x,t)=\frac{\delta(t-|x|)}{4\pi|x|};
\]
it solves the equation%
\begin{equation}
G_{tt}^{3D+}-\Delta_{x}G^{3D+}=\delta(t)\delta(|x|) \label{E:3dgreen}%
\end{equation}
subject to the radiation condition as $t\rightarrow\infty.$ We will need the
advanced Green's function $G^{3D-}(x,t),$%
\begin{equation}
G^{3D-}(x,t)=G^{3D+}(x,-t)=\frac{\delta(t+|x|)}{4\pi|x|}, \label{E:3dgreenf}%
\end{equation}
that solves (\ref{E:3dgreen}) back in time and satisfies the radiation
condition as $t\rightarrow-\infty.$

The double- and single- layer potentials are frequently used in the theory of
partial differential equations (see e.g. \cite{Vladimirov}). The wave
potentials that solve the time dependent wave equation, were introduced in
\cite{Friedman,Fulks} (see also \cite{Tuong}). We summarize below some of the
basic properties of the wave potentials that will be needed for further use.
(Here we sacrifice generality for simplicity and discuss only the properties
that will be used later in the text).

Given a plane $\Pi$ and a  function $\varphi(y,t)$, $y \in\Pi$ which is
continuous, twice differentiable in time and compactly
supported in both variables, we introduce
the (advanced)\ single layer potential $U(x,t)$ and the double-layer potential
$V(x,t)$ with the density $\varphi(y,t)$ as follows:%
\begin{align}
U(x,t)  &  =\int\limits_{\mathbb{R}^{1}}\int\limits_{\Pi}\varphi(y,t^{\prime
})G^{3D-}(y-x,t-t^{\prime})ds(y)dt^{\prime},\label{E:single3d}\\
V(x,t)  &  =\int\limits_{\mathbb{R}^{1}}\int\limits_{\Pi}\varphi(y,t^{\prime
})\frac{\partial}{\partial n_{y}}G^{3D-}(y-x,t-t^{\prime})ds(y)dt^{\prime},
\label{E:double3d}%
\end{align}
where $n$ is the normal to $\Pi$, and $ds(y)$ is the standard area element.
Note that both $U(x,t)$ and $V(x,t)$ satisfy the wave equation for all $x$
outside $\Pi$; both potentials describe waves propagating (backwards in time)
away from the support of $\varphi,$ and they satisfy the radiation condition
at infinity as $t\rightarrow-\infty.$ Also, $V(x,t)=0$ for all $x\in\Pi,$ and
due to the standard jump conditions \cite{Tuong}%
\[
\lim_{\varepsilon\rightarrow0,\varepsilon>0}V(y\pm\varepsilon n,t)=\pm\frac
{1}{2}\varphi(y,t),{\quad}y\in\Pi.
\]
Since $G^{3D-}(x,t)$ vanishes for positive values of $t,$ equation
(\ref{E:double3d}) implies that if $\varphi(y,t)$ vanishes for $t>t_{0}$ then
so does $V(x,t).$ Finally, since $n$ is constant, using (\ref{E:3dgreenf})
$V(x,t)$ can be re-written as follows%
\begin{align}
V(x,t)  &  =-\frac{1}{4\pi}n\cdot\nabla_{x}\int\limits_{\Pi}\frac
{\varphi(y,t+|y-x|)}{4\pi|y-x|}ds(y).\label{E:dontmatter}\\
&  =-\frac{1}{4\pi}\operatorname{div}\int\limits_{\Pi}n\frac{\varphi
(y,t+|y-x|)}{4\pi|y-x|}ds(y).\nonumber
\end{align}
We notice that, due to (\ref{E:dontmatter}), (for a planar surface only),
\begin{equation}
V(x,t)=-n\cdot\nabla_{x}U(x,t). \label{E:doublesingle}%
\end{equation}

In 2D, the advanced free space Green's function $G^{2D-}(x,t)$ has the
following form \cite{Vladimirov}:
\[
G^{2D-}(x,t)=\left\{
\begin{array}
[c]{cc}%
\frac{1}{2\pi\sqrt{t^{2}-x^{2}}}, & |x|<|t|,t<0,\\
0, & \mathrm{otherwise}%
\end{array}
\right.  .
\]
Further, for a given straight line $L$ and for a continuous, twice
differentiable in time, compactly supported in both variables function
$\varphi(y,t),$ $y\in L$ we define the single-layer potential $U(x,t)$ with
the density $\varphi(y,t):$%

\begin{align}
U(x,t)  &  =\int\limits_{L}\int\limits_{t}^{\infty}\varphi(y,t^{\prime}%
)G^{2D}(y-x,t-t^{\prime})dt^{\prime}dl(y)\nonumber\\
&  =\frac{1}{2\pi}\int\limits_{L}\int\limits_{|x-y|}^{\infty}\varphi
(y,t+\tau)\frac{1}{\sqrt{\tau^{2}-(x-y)^{2}}}d\tau dl(y)\nonumber\\
&  =-\frac{1}{2\pi}\int\limits_{L}\int\limits_{|x-y|}^{\infty}\left.  \left(
\frac{\partial}{\partial r}\frac{\varphi(y,r)}{r}\right)  \right\vert
_{r=t+\tau}\sqrt{\tau^{2}-(x-y)^{2}}d\tau dl(y). \label{E:single2d}%
\end{align}
Now, instead of dealing directly with $\frac{\partial}{\partial n}%
G^{2D-}(x,t)$ (which would lead to strongly singular integrals), we define the
double layer potential $V(x,t)$ by means of (\ref{E:doublesingle}) and
(\ref{E:single2d}):
\begin{align*}
V(x,t)  &  =-n\cdot\nabla_{x}U(x,t)\\
&  =\frac{1}{2\pi}n\cdot\nabla_{x}\int\limits_{L}\int\limits_{|x-y|}^{\infty
}\left.  \left(  \frac{\partial}{\partial r}\frac{\varphi(y,r)}{r}\right)
\right\vert _{r=t+\tau}\sqrt{\tau^{2}-(x-y)^{2}}d\tau dl(y)\\
&  =-\frac{1}{2\pi}\int\limits_{L}\int\limits_{|x-y|}^{\infty}\left.  \left(
\frac{\partial}{\partial r}\frac{\varphi(y,r)}{r}\right)  \right\vert
_{r=t+\tau}n\cdot\nabla_{y}\sqrt{\tau^{2}-(x-y)^{2}}d\tau dl(y)\\
&  =-\frac{1}{2\pi}\int\limits_{L}\int\limits_{|x-y|}^{\infty}\left.  \left(
\frac{\partial}{\partial r}\frac{\varphi(y,r)}{r}\right)  \right\vert
_{r=t+\tau}\frac{n\cdot(x-y)}{\sqrt{\tau^{2}-(x-y)^{2}}}d\tau dl(y),
\end{align*}
where $n$ is the normal to $L$, and $dl(y)$ is the standard arc length. As
before, $V(x,t)$ satisfies the wave equation (now in 2D) for all $x$ outside
$L$ and it describes an outgoing wave propagating (backwards in time) away
from the support of $\varphi.$ Also, $V(x,t)=0$ for all $x\in L,$ and%
\[
\lim_{\varepsilon\rightarrow0,\varepsilon>0}V(y\pm\varepsilon n,t)=\pm\frac
{1}{2}\varphi(y,t),{\quad}y\in L.
\]

\subsection{Time reversal in a cube\label{S:3drev}}

An explicit solution to the initial-boundary-value problem
(\ref{E:3dbackwards}) for the case when $\Omega$ is a cube can be obtained
using the double layer potentials, as described below.

Without loss of generality let us assume that $\Omega=(-a,a)\times
(-a,a)\times(-a,a)$, and boundary $\partial\Omega$ consists of six faces
$S_{j},$ $j=1,...,6$:
\[
\partial\Omega=\bigcup\limits_{j=1}^{6}S_{j}.
\]
In this case the wave leaves $\Omega$ after time $T=\mathrm{diam}\Omega
=2\sqrt{3}a$. Let us assume that $x=(x_{1},x_{2},x_{3}),$ and that the face
$S_{1}$ lies in the plane $x_{1}=a.$ First. we will solve the
initial-boundary-value problem with non-zero boundary conditions on face
$S_{1}$ only:%
\begin{equation}%
\begin{cases}
u_{tt}=\Delta_{x}u,\quad t\in\lbrack0,T],\quad x\in\Omega\subset\mathbb{R}%
^{3},\\
u(x,T)=0,\quad u_{t}(x,T)=0,\\
\lim_{\varepsilon\rightarrow0,\varepsilon>0}u(y-\varepsilon n_{1}%
,t)=P(y,t),\quad y\in S_{1},\\
u(y,t)=0,\quad y\in\bigcup\limits_{j=2}^{6}S_{j},
\end{cases}
\label{E:oneface3d}%
\end{equation}
where $n_{1}$ is the exterior normal to $S_{1}.$ The solution procedure we are
about to outline will apply as well to the other faces $S_{j},$ $j=2,...,6$,
and the solution of problem (\ref{E:3dbackwards}) will then be readily
obtained as the sum of the solutions corresponding to each face.

To simplify the presentation let us introduce an operator $E^{2D}$ that
extends any function $h(x_{2},x_{3})$ defined in the open square
$(-a,a)\times(-a,a)$ to a function $h^{\mathrm{repl}}(x_{2},x_{3})$ defined in
$\mathbb{R}^{2}$ by means of odd reflections. In detail, if $h^{\mathrm{repl}%
}=E^{2D}h,$ then $h^{\mathrm{repl}}$ coincides with $h$ within $(-a,a)\times
(-a,a);$ further, for any point $(x_{2},x_{3})$ lying in an open square
$s_{n,m}$ with the side $2a$ and centered at $(2na,2ma),$ $h^{\mathrm{repl}%
}(x_{2},x_{3})$ is defined by the following formula:
\[
h^{\mathrm{repl}}(x_{2},x_{3})=(-1)^{m+n}h(a,(-1)^{n}(x_{2}-2na),(-1)^{m}%
(x_{3}-2ma)).
\]
Finally, on the boundary of each square $s_{n,m}$ we define function
$h^{\mathrm{repl}}(x_{2},x_{3})$ to be zero.

Let us denote by $\Pi_{1}^{a}$ the plane containing $S_{1}$, and introduce
plane $\Pi_{1}^{b}$ parallel to $\Pi_{1}^{a}$ and passing through the point
$(-3a,0,0).$ We equip both planes with the same normal $n_{1}$ coinciding with
the exterior normal to face $S_{1}.$ We introduce function $P_{1}%
^{\mathrm{repl}}(x_{2},x_{3},t)=E^{2D}\left(  \left.  P(x,t)\right\vert
_{S_{1}}\right)  ,$ i.e. we replicate the boundary values of $u$ on face
$S_{1}.$ Now, for $t\in\lbrack0,T],$ we define the double-layer potentials
$V_{1}^{a}(x,t)$ and $V_{1}^{b}(x,t)$ with the density $-2P_{1}^{\mathrm{repl}%
}(y_{2},y_{3},t)$ supported on planes $\Pi_{1}^{a}$ and $\Pi_{1}^{b}$:%
\begin{equation}%
\begin{array}
[c]{c}%
V_{1}^{a}(x,t)=\frac{1}{2\pi}\operatorname{div}\int\limits_{\mathbb{R}^{2}%
}\frac{n_{1}}{|(a,y_{2},y_{3})-x|}P_{1}^{\mathrm{repl}}(y_{2},y_{3}%
,t+|(a,y_{2},y_{3})-x|)dy_{2}dy_{3},\\
V_{1}^{b}(x,t)=\frac{1}{2\pi}\operatorname{div}\int\limits_{\mathbb{R}^{2}%
}\frac{n_{1}}{|(-3a,y_{2},y_{3})-x|}P_{1}^{\mathrm{repl}}(y_{2},y_{3}%
,t+|(-3a,y_{2},y_{3})-x|)dy_{2}dy_{3}.
\end{array}
\label{E:doublelayers}%
\end{equation}
Note that since $P_{1}^{\mathrm{repl}}$ is finitely supported in $[0,T]$ in
the time variable, the integrands in the above expressions vanish for large
values of $y_{2}$ or $y_{3},$ and therefore the integration is actually
performed over bounded subsets of $\mathbb{R}^{2}$ lying within the ball
$B(x,T-t).$ In order to make expressions (\ref{E:doublelayers}) simpler, we
will abuse notation by writing $P_{1}^{\mathrm{repl}}(y,t)$ for $y\in\Pi
_{1}^{a}$ or $y\in\Pi_{1}^{b}$ instead of $P_{1}^{\mathrm{repl}}(y_{2}%
,y_{3},t):$%
\begin{equation}
V_{1}^{a,b}(x,t)=\frac{1}{2\pi}\operatorname{div}\int\limits_{\Pi_{1}%
^{a,b}\cap B(x,T-t)}\frac{n_{1}}{|y-x|}P_{1}^{\mathrm{repl}}(y,t+|y-x|)ds(y),
\label{E:doublelayers3d}%
\end{equation}
where $ds(y)$ is the standard area element, and $V_{1}^{a,b}$ means either
$V_{1}^{a}$ or $V_{1}^{b}.$

The combined potential $V_{1}(x,t)=V_{1}^{a}(x,t)+V_{1}^{b}(x,t)$ has some
interesting properties. Let us denote by $S_{2}$ the face of cube $\Omega$
opposite to $S_{1}.$ We notice that planes $\Pi_{1}^{a}$ and $\Pi_{1}^{b}$ are
symmetric with respect to $S_{2}$, and, due to the choice of the normals,
potential $V_{1}(x,t)$ vanishes on $S_{2}.$ Next, due to the oddness of
$P_{1}^{\mathrm{repl}},$ potential $V_{1}(x,t)$ vanishes on the other four
faces $S_{j}$ of the cube $(j=3,...,6$). Finally, since $T$ is smaller than
the distance between $\Pi_{1}^{a}$ and $\Pi_{1}^{b}$, potential $V_{1}%
^{b}(x,t)$ equals 0 on $\Pi_{1}^{a}$ for all $t\leq T.$ Therefore,
\[
\lim_{\varepsilon\rightarrow0,\varepsilon>0}V_{1}(y-\varepsilon n_{1}%
,t)=\lim_{\varepsilon\rightarrow0,\varepsilon>0}V_{1}^{a}(y-\varepsilon
n_{1},t)=P(y,t),\quad y\in S_{1}.
\]
Thus, potential $V_{1}(x,t)$ solves initial-boundary-value problem
(\ref{E:oneface3d}). For brevity, we will re-write the expression for
$V_{1}(x,t)$ in the following form%
\[
V_{1}(x,t)=\frac{1}{2\pi}\operatorname{div}\int\limits_{(\Pi_{1}^{a}\cup
\Pi_{1}^{b})\cap B(x,T-t)}\frac{n_{1}}{|y-x|}P_{1}^{\mathrm{repl}%
}(y,t+|y-x|)ds(y).
\]

In order to obtain the desired solution to (\ref{E:3dbackwards}), we replicate
the boundary values on the faces $S_{2},...,S_{6}:$%
\begin{equation}
P_{j}^{\mathrm{repl}}(\cdot,\cdot,t)=E^{2D}\left(  \left.  P(\cdot
,t)\right\vert _{S_{j}}\right)  ,j=2,...,6. \label{E:repldata3d}%
\end{equation}
Next, we introduce planes $\Pi_{j}^{a},$ $\Pi_{j}^{b}$ with normals $n_{j},$
and define double layer potentials $V_{j}$ with the densities $-2P_{j}%
^{\mathrm{repl}}(y,t)$ on $\Pi_{j}^{a},$ $\Pi_{j}^{b}$:%
\[
V_{j}(x,t)=\frac{1}{2\pi}\operatorname{div}\int\limits_{(\Pi_{j}^{a}\cup
\Pi_{j}^{b})\cap B(x,T-t)}\frac{n_{j}}{|y-x|}P_{j}^{\mathrm{repl}%
}(y,t+|y-x|)ds(y),\qquad j=2,...,6.
\]
The sum of these potentials $V(x,t)=$ $\sum\nolimits_{j=1}^{6}V_{j}(x,t)$
solves the initial-boundary-value problem (\ref{E:3dbackwards}), and,
therefore, $f(x)$ equals to $V(x,0):$%
\begin{equation}
f(x)=\sum\limits_{j=1}^{6}V_{j}(x,0)=\frac{1}{2\pi}\operatorname{div}%
\sum\limits_{j=1}^{6}\int\limits_{(\Pi_{j}^{a}\cup\Pi_{j}^{b})\cap
B(x,T)}n_{j}\frac{P_{j}^{\mathrm{repl}}(y,|y-x|)}{|y-x|}ds(y).
\label{E:3dwaveinversion}%
\end{equation}
Interestingly, while each of the potentials $V_{j}^{a,b}(x,t)$ solves the wave
equation in the sense of distributions (due to the discontinuities in the
edges of $\Omega)$ the sum $\sum\limits_{j=1}^{6}V_{j}(x,0)$ of these
potentials solves the wave equation in the classical sense (if $f(x)\in
C_{0}^{2}(\Omega)).$ Thus, we have proven the following

\begin{theorem}
The $C_{0}^{2}(\Omega)$ initial condition $f(x)$ of the initial-boundary-value
problem (\ref{E:wave_eq}) in the case $n=3$ can be reconstructed from the
boundary values $P(y,t)$ (defined by (\ref{E:wavedata})) by formulae
(\ref{E:3dwaveinversion}) and (\ref{E:repldata3d}).
\end{theorem}

Formula (\ref{E:3dwaveinversion}) closely resembles the so-called "universal
backprojection formula" \cite{XuWang} especially the version applicable to the
measurements done over an infinite plane. However, in the present case the
measurements are performed over a bounded surface (that of the cube), and the
exact reconstruction at any point $x$ is obtained by integration over a
bounded set $B(x,T)\cap\bigcup\nolimits_{j=1}^{6}(\Pi_{j}^{a}\cup\Pi_{j}%
^{b}).$

\subsection{Time reversal in a square\label{S:2drev}}

A similar technique can be used to find the explicit solution to the
initial-boundary-value problem (\ref{E:back2d}) in the case when $\Omega$ is a
square. The geometry of the problem is actually easier, thanks to the lower
dimensionality of the problem. However, in 2D the Huygens principle does not
hold and the solution will result from the integration over unbounded sets.

We will assume that $\Omega$ is a square $(-a,a)\times(-a,a)$, $S_{1}$ is the
side of $\Omega$ corresponding to $x_{1}=a,$ and $S_{2}$ is the side contained
by the line $x_{1}=-a.$

First, let us find solution $u(x,t)$ to the following initial-boundary-value problem:%

\begin{equation}%
\begin{cases}
u_{tt}=\Delta_{x}u,\quad t\in\lbrack0,T],\quad x\in\Omega\subset\mathbb{R}%
^{2},\\
u(x,T)=0,\quad u_{t}(x,T)=0,\\
\lim\limits_{\varepsilon\rightarrow0,\varepsilon>0}u(y-\varepsilon
n,t)=P(y,t)\eta(T,t),\quad y\in S_{1},\\
u(y,t)=0,\quad y\in S_{j},\quad j=2,3,4,
\end{cases}
\label{E:truncated2d}%
\end{equation}
where $\eta(T,t)$ is a $C^{\infty}$ cut-off function equal 1 for $t\in
\lbrack0,T-1]$ and vanishing with all the derivatives at $t=T.$ We introduce
operator $E^{1D}$ that extends a function $h(x)$ defined in the interval
$(-a,a),$ to a function $h^{\mathrm{repl}}(x)$ defined on $\mathbb{R}^{1}$ by
means of odd reflections. In detail, within the interval $(-a,a)$ function
$h^{\mathrm{repl}}$ coincides with $h.$\ For any $x\neq(2m+1)a,$
$m\in\mathbb{Z},$ there is a unique integer $n(x)$ such that $|x-2na|<a.$ We
define $h^{\mathrm{repl}}(x)=(-1)^{n(x)}h((-1)^{n(x)}(x-2n(x)a)).$ Finally, at
the odd integer points $h^{\mathrm{repl}}$ vanishes: $h^{\mathrm{repl}%
}((2m+1)a)=0,$ $\forall m\in\mathbb{Z}.$

Let us call $L_{1}$ the straight line containing $S_{1}.$ Introduce the family
of lines $\Phi_{1}$ consisting of all straight lines parallel to $L_{1}$, and
passing through the points $(a+4ka,0),$ $k\in\mathbb{Z}$. (Obviously,
$L_{1}\in\Phi_{1}).$ The exterior normal $n_{1}$ to the side $S_{1}$ of the
square will be used as the normal to all the lines in $\Phi_{1}$. Let us
define function $P_{1}^{\mathrm{repl}}(x_{2},t)$ by replicating values of
$P(x,t)$ for $x\in S_{1}:$%
\[
P_{1}^{\mathrm{repl}}(x_{2},t)=E^{1D}\left(  P((a,x_{2}),t)\right)  .
\]
We will abuse notation by writing $P_{1}^{\mathrm{repl}}(x,t)$ instead of the
more accurate $P_{1}^{\mathrm{repl}}(x_{2}(x),t).$ Further, let us introduce
the double layer potential $V_{1}(x,t,T)$ with the density $-2P_{1}%
^{\mathrm{repl}}(x_{2},t)\eta(T,t)$ supported on all lines in the family
$\Phi_{1},$ defined by the formula%
\[
V_{1}(x,t,T)=\frac{1}{\pi}\int\limits_{\Phi_{1}}\int\limits_{|x-y|}^{\infty
}\left.  \left(  \frac{\partial}{\partial r}\frac{P_{1}^{\mathrm{repl}%
}(y,r)\eta(T,r)}{r}\right)  \right\vert _{r=t+\tau}\frac{n_{1}\cdot
(x-y)}{\sqrt{\tau^{2}-(x-y)^{2}}}d\tau dl(y).
\]
Note that in the above formula the integration in $y$ is actually performed
over a bounded subset of $\Phi_{1}$ due to the finite support of $\eta$ in the
second variable.

Potential $V_{1}(x,t,T)$ satisfies the wave equation in $\Omega.$ Since
$P_{1}^{\mathrm{repl}}$ is odd with respect to the straight lines containing
all sides of $\partial\Omega,$ $V_{1}^{T}(x,t,T)$ vanishes on $\partial
\Omega.$ However, $V_{1}(x,t,T)$ is not continuous across $S_{1}.$ The jump is
such that%
\[
\lim_{\varepsilon\rightarrow0,\varepsilon>0}V_{1}(y-\varepsilon n_{1}%
,t,T)=P_{1}^{\mathrm{repl}}(y,t)\eta(T,t)=P(y,t)\eta(T,t),\quad y\in
S_{1},t\in\lbrack0.T].
\]
Therefore, $V_{1}(x,t,T)$ solves the problem (\ref{E:truncated2d}).

The next step is to obtain a solution $u^{T}(x,t)$ of the following boundary
problem:%
\begin{equation}%
\begin{cases}
u_{tt}^{T}=\Delta_{x}u^{T},\quad t\in\lbrack0,T],\quad x\in\Omega
\subset\mathbb{R}^{2},\\
u^{T}(x,T)=0,\quad u_{t}^{T}(x,T)=0,\\
u^{T}(y,t)=P(y,t)\eta(T,t),\quad y\in\partial\Omega.
\end{cases}
\label{E:truncated2dall}%
\end{equation}
This can be achieved by defining functions $P_{j}^{\mathrm{repl}}%
(\cdot,t)=E^{1D}\left(  \left.  P(\cdot,t)\right\vert _{S_{j}}\right)  $ and
families of parallel lines $\Phi_{j}$ with normals $n_{j},$ $j=2,3,4,$ and by
introducing the corresponding double layer potentials $V_{j}(x,t,T),$ in a
fashion similar to our definition of $P_{1}^{\mathrm{repl}},\Phi_{1},n_{1},$
and $V_{1}(x,t,T).$ Now each $V_{j}(x,t,T)$ satisfies the boundary problem
with the inhomogeneous conditions on the side $j,$ and the homogeneous
conditions on the other three sides of the square. In turn, the sum of these
potentials $u^{T}$%
\begin{align*}
u^{T}(x,t)  &  =\sum_{j=1}^{4}V_{j}(x,t,T)=\frac{1}{\pi}\operatorname{div}%
\sum_{j=1}^{4}\left(  \int\limits_{0}^{T-t}\int\limits_{\underset{|y-x|<\tau
}{\Phi_{j}}}n_{j}\frac{P_{j}^{\mathrm{repl}}(y,t+\tau)\eta(T,t+\tau)}%
{\sqrt{\tau^{2}-(y-x)^{2}}}dl(y)d\tau\right) \\
&  =\frac{1}{\pi}\sum_{j=1}^{4}\int\limits_{\Phi_{j}}\int\limits_{|x-y|}%
^{\infty}\left.  \left(  \frac{\partial}{\partial r}\frac{P_{j}^{\mathrm{repl}%
}(y,r)\eta(T,r)}{r}\right)  \right\vert _{r=t+\tau}\frac{n_{j}\cdot
(x-y)}{\sqrt{\tau^{2}-(x-y)^{2}}}d\tau dl(y)
\end{align*}
solves (\ref{E:truncated2dall}).

It is known \cite{Hristova} that, as $T$ grows to infinity, $u^{T}(x,0)$
converges uniformly to $f(x)$:%
\begin{align}
f(x)  &  =\lim_{T\rightarrow\infty}u^{T}(x,0)=\lim_{T\rightarrow\infty}%
\frac{1}{\pi}\sum_{j=1}^{4}\int\limits_{\Phi_{j}}\int\limits_{|x-y|}^{\infty
}\left(  \frac{\partial}{\partial r}\frac{P_{j}^{\mathrm{repl}}(y,r)\eta
(T,r)}{r}\right)  \frac{n_{j}\cdot(x-y)}{\sqrt{r^{2}-(x-y)^{2}}}%
drdl(y)\label{E:almostfinalwave2D}\\
&  =\frac{1}{\pi}\sum_{j=1}^{4}\int\limits_{\Phi_{j}}\int\limits_{|x-y|}%
^{\infty}\left(  \frac{\partial}{\partial r}\frac{P_{j}^{\mathrm{repl}}%
(y,r)}{r}\right)  \frac{n_{j}\cdot(x-y)}{\sqrt{r^{2}-(x-y)^{2}}}drdl(y)
\label{E:finalwave2D}%
\end{align}
The latter formula is very similar to the explicit inversion formula for
measurements done from an infinite line \cite{Burgh}. Formula
(\ref{E:finalwave2D}) can also be re-written using the Green's function
$G^{2D-}(x-y,r)$ to illustrate its connection to the surface potentials:%
\[
f(x)=\frac{1}{\pi}\sum_{j=1}^{4}\int\limits_{\Phi_{j}}\int\limits_{0}^{\infty
}n_{j}\cdot(x-y)\left(  \frac{\partial}{\partial r}\frac{P_{j}^{\mathrm{repl}%
}(y,r)}{r}\right)  G^{2D-}(x-y,r)drdl(y).
\]

\begin{theorem}
The $C_{0}^{2}(\Omega)$ initial condition $f(x)$ of the initial-boundary-value
problem (\ref{E:wave_eq}) in the case $n=2$ can be reconstructed from the
boundary values $P(y,t)$ (defined by (\ref{E:wavedata})) by formulae
(\ref{E:finalwave2D}).
\end{theorem}

\section{Explicit inversion of the spherical means: a cube and a
square\label{S:means}}

\subsection{Inversion formula in 3D\label{S:3dmeans}}

In 3D, in order to reconstruct $f(x)$ from its spherical means $M(y,t)$ with
centers on $\partial\Omega,$ we introduce replicated spherical means
$M_{j}^{\mathrm{repl}}(y,t)$ defined as follows%
\begin{equation}
M_{j}^{\mathrm{repl}}(\cdot,\cdot,t)=E^{2D}\left(  \left.  M(\cdot
,t)\right\vert _{S_{j}}\right)  ,j=1,...,6. \label{E:replaverdef}%
\end{equation}
Since the replication operator $E^{2D}$ is invariant with respect to the time
variable, from (\ref{E:3daverages}) we conclude that%
\begin{equation}
{P}_{j}^{\mathrm{repl}}{(y,t)=}\frac{\partial}{\partial t}tM_{j}%
^{\mathrm{repl}}(y,t). \label{E:replaver}%
\end{equation}
Now formula for reconstructing $f(x)$ from spherical averages in 3D is
obtained by substituting (\ref{E:replaver}) into (\ref{E:3dwaveinversion}):%
\begin{equation}
f(x)=\frac{1}{2\pi}\operatorname{div}\sum\limits_{j=1}^{6}\int\limits_{(\Pi
_{j}^{a}\cup\Pi_{j}^{b})\cap B(x,T)}n_{j}\left.  \left(  \frac{1}{t}%
\frac{\partial}{\partial t}tM_{j}^{\mathrm{repl}}(y,t)\right)  \right\vert
_{t=|x-y|}ds(y). \label{E:3dmeaninversion}%
\end{equation}
We thus have proven

\begin{theorem}
In 3D, a $C_{0}^{2}(\Omega)$ function $f(x)$ supported within the cube
$\Omega$ can be reconstructed from its spherical means $M(y,t)$ with centers
on $\partial\Omega$ by the formula (\ref{E:3dmeaninversion}), where
$M_{j}^{\mathrm{repl}}(y,t)$ is defined by (\ref{E:replaverdef}).
\end{theorem}

Formula (\ref{E:3dmeaninversion}) has a form of the filtration by
differentiation, followed by the backprojection, followed by the divergence
operator. It is similar to the inversion formula obtained by the author
\cite{kunya} for the inversion of the spherical means with centers on a sphere
(in 3D). In the latter case the integration is over the surface of a sphere,
and no replication of the data is needed.

\subsection{Inversion formula in 2D\label{S:2dmeans}}

In this section we derive an explicit inversion formula for finding $f(x)$
from the values of its circular means with the centers lying on the boundary
of a square. The starting point for the derivation is equation
(\ref{E:finalwave2D}) that reconstructs the sought function from the boundary
values of the solution of the wave equation; we re-write it in the following form%

\begin{align}
f(x)  &  =\frac{1}{\pi}\sum_{j=1}^{4}\int\limits_{\Phi_{j}}n_{j}%
\cdot(x-y)\left.  \left[  \int\limits_{s}^{\infty}\left(  \frac{\partial
}{\partial r}\frac{P_{j}^{\mathrm{repl}}(y,r)}{r}\right)  \frac{1}{\sqrt
{r^{2}-s^{2}}}dr\right]  \right\vert _{s=|x-y|}dl(y)\nonumber\\
&  =\frac{1}{\pi}\sum_{j=1}^{4}\int\limits_{\Phi_{j}}n_{j}\cdot(x-y)\hat
{P}_{j}^{\mathrm{repl}}(y,|x-y|)dl(y) \label{E:2deq2}%
\end{align}
where%
\begin{equation}
\hat{P}_{j}^{\mathrm{repl}}(y,s)=\int\limits_{s}^{\infty}\left(
\frac{\partial}{\partial r}\frac{P_{j}^{\mathrm{repl}}(y,r)}{r}\right)
\frac{1}{\sqrt{r^{2}-s^{2}}}dr. \label{E:2dder}%
\end{equation}
First, we will express $\hat{P}_{j}^{\mathrm{repl}}(y,s)$ in terms of the
circular averages $m(y,r).$ Let us apply the extension operator to $m(y,r)$
and define $m_{j}^{\mathrm{repl}}(y,r)$ as follows:%
\begin{equation}
m_{j}^{\mathrm{repl}}(\cdot,r)=E^{1D}\left(  \left.  m(\cdot,r)\right\vert
_{S_{j}}\right)  ,\qquad j=1,...,4. \label{E:2drepl}%
\end{equation}
Function $P_{j}^{\mathrm{repl}}(\cdot,t)$ is related to $m_{j}^{\mathrm{repl}%
}(\cdot,r)$ in the same way as $P(\cdot,t)$ is related to $m(\cdot,r)$:%
\begin{align*}
P_{j}^{\mathrm{repl}}(y,t)  &  =\frac{\partial}{\partial t}\int\limits_{0}%
^{t}\frac{m_{j}^{\mathrm{repl}}(y,r)rdr}{\sqrt{t^{2}-r^{2}}}=\frac{\partial
}{\partial t}\int\limits_{0}^{t}\sqrt{t^{2}-r^{2}}\frac{\partial}{\partial
r}m_{j}^{\mathrm{repl}}(y,r)dr=\int\limits_{0}^{t}\frac{t\frac{\partial
}{\partial r}m_{j}^{\mathrm{repl}}(y,r)}{\sqrt{t^{2}-r^{2}}}dr\\
&  =\int\limits_{0}^{t}\frac{t\frac{\partial}{r\partial r}m_{j}^{\mathrm{repl}%
}(y,r)}{\sqrt{t^{2}-r^{2}}}rdr=\int\limits_{0}^{t}\frac{\partial}{\partial
r}\left(  \frac{\partial}{r\partial r}m_{j}^{\mathrm{repl}}(y,r)\right)
t\sqrt{t^{2}-r^{2}}dr,
\end{align*}
where we integrated by parts twice, and used the fact that $m_{j}%
^{\mathrm{repl}}(y,0)=0$ for all $y.$ Now one can compute the derivative
needed in (\ref{E:2dder}):%
\[
\frac{\partial}{\partial t}\left(  P_{j}^{\mathrm{repl}}(y,t)/t\right)
=\frac{\partial}{\partial t}\int\limits_{0}^{t}\frac{\partial}{\partial
r}\left(  \frac{\partial}{r\partial r}m_{j}^{\mathrm{repl}}(y,r)\right)
\sqrt{t^{2}-r^{2}}dr=t\int\limits_{0}^{t}\frac{\frac{\partial}{\partial
r}\left(  \frac{\partial}{r\partial r}m_{j}^{\mathrm{repl}}(y,r)\right)
}{\sqrt{t^{2}-r^{2}}}dr.
\]
By substituting the above equation into (\ref{E:2dder}), interchanging the
integration order, and noticing that $m_{j}^{\mathrm{repl}}(y,r)=0$ for
$r>2\sqrt{2}a,$ one obtains the following expression for $\hat{P}%
_{j,T}^{\mathrm{repl}}(y,s):$%
\begin{align}
\hat{P}_{j,T}^{\mathrm{repl}}(y,s)  &  =\lim_{T\rightarrow\infty}%
\int\limits_{s}^{T}\int\limits_{0}^{t}\frac{\partial}{\partial r}\left(
\frac{\partial}{r\partial r}m_{j}^{\mathrm{repl}}(y,r)\right)  \frac{t}%
{\sqrt{t^{2}-r^{2}}\sqrt{t^{2}-s^{2}}}drdt\nonumber\\
&  =\lim_{T\rightarrow\infty}\int\limits_{0}^{2\sqrt{2}a}\frac{\partial
}{\partial r}\left(  \frac{\partial}{r\partial r}m_{j}^{\mathrm{repl}%
}(y,r)\right)  \left[  \int\limits_{\max(s,r)}^{T}\frac{t}{\sqrt{t^{2}-r^{2}%
}\sqrt{t^{2}-s^{2}}}dt\right]  dr\nonumber\\
&  =\lim_{T\rightarrow\infty}\int\limits_{0}^{2\sqrt{2}a}\frac{\partial
}{\partial r}\left(  \frac{\partial}{r\partial r}m_{j}^{\mathrm{repl}%
}(y,r)\right)  \left[  \ln\left(  \sqrt{T^{2}-s^{2}}+\sqrt{T^{2}-r^{2}%
}\right)  -\frac{1}{2}\ln\left(  |r^{2}-s^{2}|\right)  \right]  dr\nonumber\\
&  =-\frac{1}{2}\int\limits_{0}^{2\sqrt{2}a}\ln\left(  |r^{2}-s^{2}|\right)
\frac{\partial}{\partial r}\left(  \frac{\partial}{r\partial r}m_{j}%
^{\mathrm{repl}}(y,r)\right)  dr\nonumber\\
&  +\lim_{T\rightarrow\infty}\int\limits_{0}^{2\sqrt{2}a}\frac{\partial
}{\partial r}\left(  \frac{\partial}{r\partial r}m_{j}^{\mathrm{repl}%
}(y,r)\right)  \ln\left(  \sqrt{T^{2}-s^{2}}+\sqrt{T^{2}-r^{2}}\right)  dr
\label{E:2deq1}%
\end{align}
The last integral in (\ref{E:2deq1}) vanishes:%
\begin{align*}
&  \lim_{T\rightarrow\infty}\int\limits_{0}^{2\sqrt{2}a}\frac{\partial
}{\partial r}\left(  \frac{\partial}{r\partial r}m_{j}^{\mathrm{repl}%
}(y,r)\right)  \ln\left(  \sqrt{T^{2}-s^{2}}+\sqrt{T^{2}-r^{2}}\right)  dr\\
&  =\lim_{T\rightarrow\infty}\int\limits_{0}^{2\sqrt{2}a}\frac{\partial
}{r\partial r}m_{j}^{\mathrm{repl}}(y,r)\frac{r}{\left(  \sqrt{T^{2}-s^{2}%
}+\sqrt{T^{2}-r^{2}}\right)  \sqrt{T^{2}-r^{2}}}dr=0,
\end{align*}
where we integrated by parts and used the finite support of $m_{j}%
^{\mathrm{repl}}(y,r)$ in $r.$ Now (\ref{E:2deq1}) can be re-written in the
following form%
\[
\hat{P}_{j,T}^{\mathrm{repl}}(y,s)=-\frac{1}{2}\int\limits_{0}^{2\sqrt{2}a}%
\ln\left(  |r^{2}-s^{2}|\right)  \frac{\partial}{\partial r}\left(
\frac{\partial}{r\partial r}m_{j}^{\mathrm{repl}}(y,r)\right)  dr=\mathrm{PV}%
\int\limits_{0}^{2\sqrt{2}a}\frac{1}{r^{2}-s^{2}}\frac{\partial}{\partial
r}m_{j}^{\mathrm{repl}}(y,r)dr,
\]
where "PV" indicates that the integral is understood in the principal value
sense. Substitution of the latter expression for $\hat{P}_{j,T}^{\mathrm{repl}%
}(y,s)$ into (\ref{E:2deq2}) yields the following sequence of equivalent
reconstruction formulae for $f(x)$:%

\begin{align}
f(x)  &  =\frac{1}{\pi}\sum_{j=1}^{4}\int\limits_{\Phi_{j}}n_{j}%
\cdot(x-y)\left.  \left(  \mathrm{PV}\int\limits_{0}^{2\sqrt{2}a}\frac
{1}{r^{2}-s^{2}}\frac{\partial}{\partial r}m_{j}^{\mathrm{repl}}%
(y,r)dr\right)  \right\vert _{s=|x-y|}dl(y)\nonumber\\
&  =-\frac{1}{2\pi}\sum_{j=1}^{4}\int\limits_{\Phi_{j}}\frac{n_{j}\cdot
(x-y)}{|x-y|}\left.  \left(  \frac{\partial}{\partial s}\int\limits_{0}%
^{2\sqrt{2}a}\ln(|r^{2}-s^{2}|)\frac{\partial}{\partial r}m_{j}^{\mathrm{repl}%
}(y,r)dr\right)  \right\vert _{s=|x-y|}dl(y)\nonumber\\
&  =\frac{1}{\pi}\sum_{j=1}^{4}\int\limits_{\Phi_{j}}\frac{n_{j}\cdot
(x-y)}{|x-y|}\left.  \left(  \frac{\partial}{\partial s}\mathrm{PV}%
\int\limits_{0}^{2\sqrt{2}a}\frac{rm_{j}^{\mathrm{repl}}(y,r)}{r^{2}-s^{2}%
}dr\right)  \right\vert _{s=|x-y|}dl(y)\label{E:2dmean1}\\
&  =\frac{1}{\pi}\sum_{j=1}^{4}\int\limits_{\Phi_{j}}n_{j}\cdot\nabla
_{x}\left.  \left(  \mathrm{PV}\int\limits_{0}^{2\sqrt{2}a}\frac
{rm_{j}^{\mathrm{repl}}(y,r)}{r^{2}-s^{2}}dr\right)  \right\vert
_{s=|x-y|}dl(y)\label{E:2dmean2}\\
&  =\frac{1}{\pi}\operatorname{div}\sum_{j=1}^{4}\int\limits_{\Phi_{j}}%
n_{j}\left[  \mathrm{PV}\int\limits_{0}^{2\sqrt{2}a}\frac{m_{j}^{\mathrm{repl}%
}(y,r)}{r^{2}-(x-y)^{2}}rdr\right]  dl(y). \label{E:2dmean3}%
\end{align}
Thus, the following statement holds:

\begin{theorem}
In 2D, a $C_{0}^{2}(\Omega)$ function $f(x)$ supported within the square
$\Omega$ can be reconstructed from its circular means $m(y,r)$ with centers on
$\partial\Omega$ by the formula (\ref{E:2dmean3}) (or, equivalently, by
formulae (\ref{E:2dmean1}) or (\ref{E:2dmean2})), where $m_{j}^{\mathrm{repl}%
}(y,r)$ is defined by (\ref{E:2drepl}).
\end{theorem}


Formula (\ref{E:2dmean3}) closely resembles formula 8.19 \cite{AKK} that
reconstructs a function from its circular means centered on a circle.

The drawback of formulae (\ref{E:2dmean1})-(\ref{E:2dmean3}) is that the
integration has to be done over an unbounded set $\cup_{j=1}^{4}\Phi_{j}.$
However, by analyzing the expression in parentheses in (\ref{E:2dmean1}) one
can notice that the integrand of the outer integral decreases as $|x-y|^{-3}$
for large values of $y.$ Moreover, for values of $|x-y|>2\sqrt{2}a$ the
integrand is infinitely smooth. Therefore, by integrating over a subset
$\Phi^{\mathrm{trunc}}$  containing all points $y\in\cup_{j=1}^{4}\Phi_{j}$
such that $\mathrm{dist}(y,\Omega)\leq\mathrm{diam}\Omega$ one correctly
reconstructs the singularities of $f(x),$ and obtains a good quantitative
approximation to $f(x).$ This conclusion is supported by the numerical
evidence presented in the next section.

\subsection{Reconstruction in the presence of exterior
sources\label{S:exterior}}

The time reversal reconstruction methods in TAT/PAT (i.e. methods based on the
numerical solution of the wave equation backwards in time) have the following
interesting property. If the support of the initial perturbation $f(x)$ has a
part lying outside the region enclosed by the acquisition surface
$\partial\Omega$, the time reversal still correctly reconstructs the part of
$f(x)$ lying within $\Omega$. In other words, the exterior sources do not
affect the reconstruction within $\Omega$ \cite{kunyaser,AK}. Since
reconstruction methods based on the expansions in the series of the
eigenfunctions of the Dirichlet Laplacian on $\Omega$ \cite{kunyaser,AK} are
theoretically equivalent to the time reversal, they also inherit the
insensitivity to the exterior sources. On the other hand, all previously known
filtration/backprojection-type inversion formulae for the spherical mean
transform do not have this property: in the presence of an exterior source the
reconstruction within $\Omega$\ will be incorrect.

The present inversion formulae (and their generalizations in the following
sections) are based on the explicit solution of the wave equation (by means of
double layer potentials). Thus, they inherit from the time reversal methods
the insensitivity to the exterior sources. A numerical example demonstrating
accurate reconstruction in the presence of an exterior source is presented in
the next section.

\subsection{Numerical examples\label{S:numerics}}

\begin{figure}[t]
\begin{center}
\subfigure[]{\includegraphics[width=1.8in,height=1.8in]{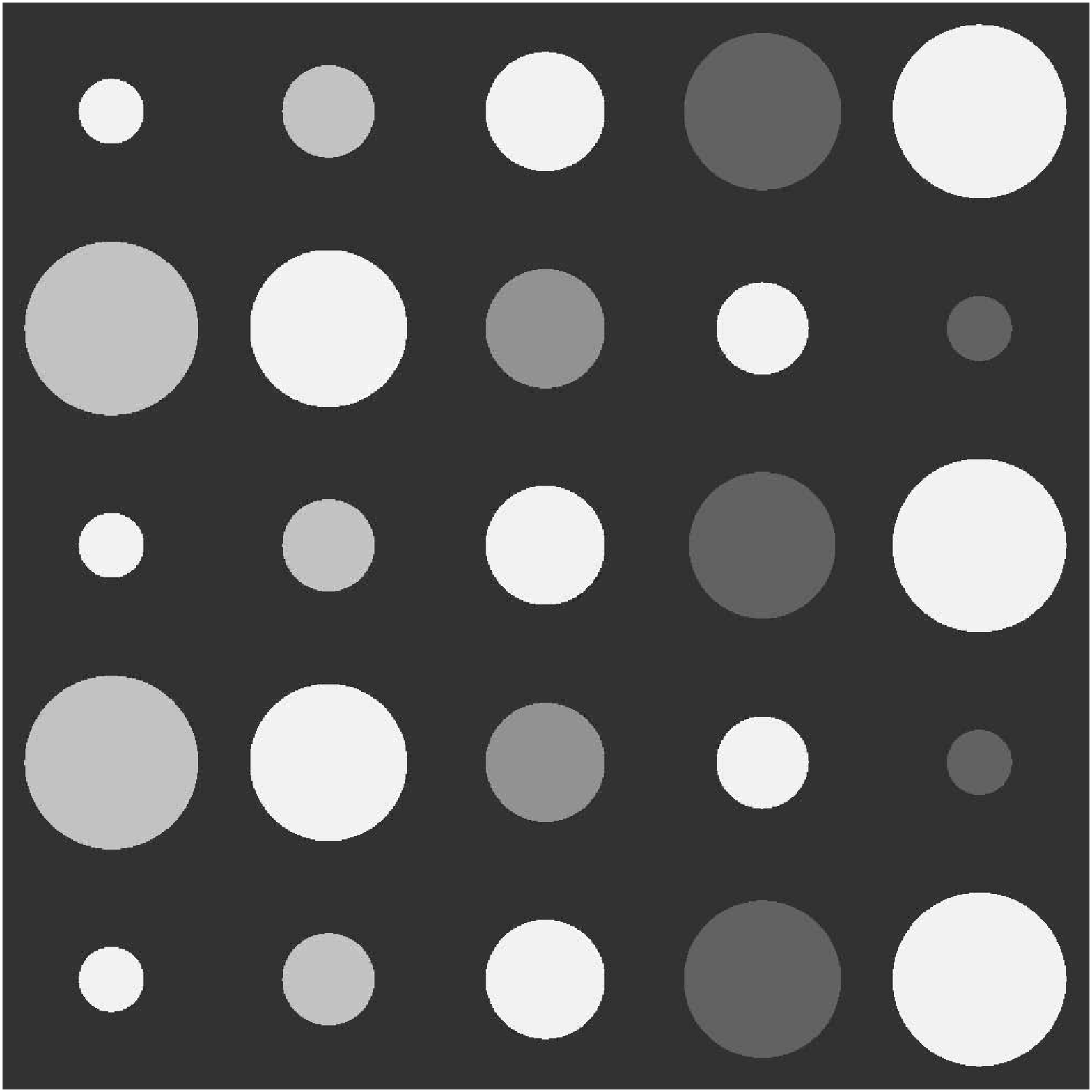} }
\subfigure[]{\includegraphics[width=1.8in,height=1.8in]{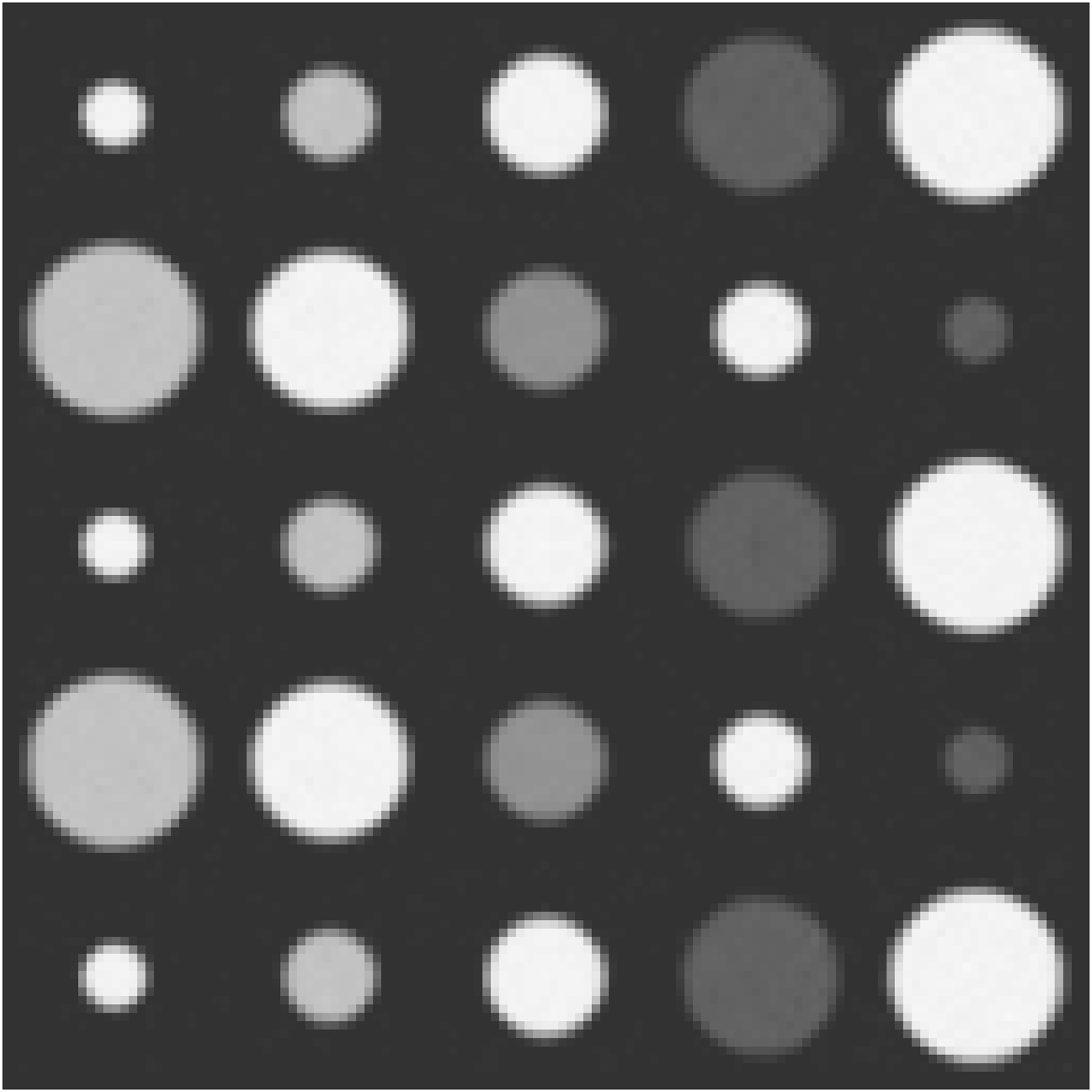}      }
\subfigure[]{\includegraphics[width=1.8in,height=1.8in]{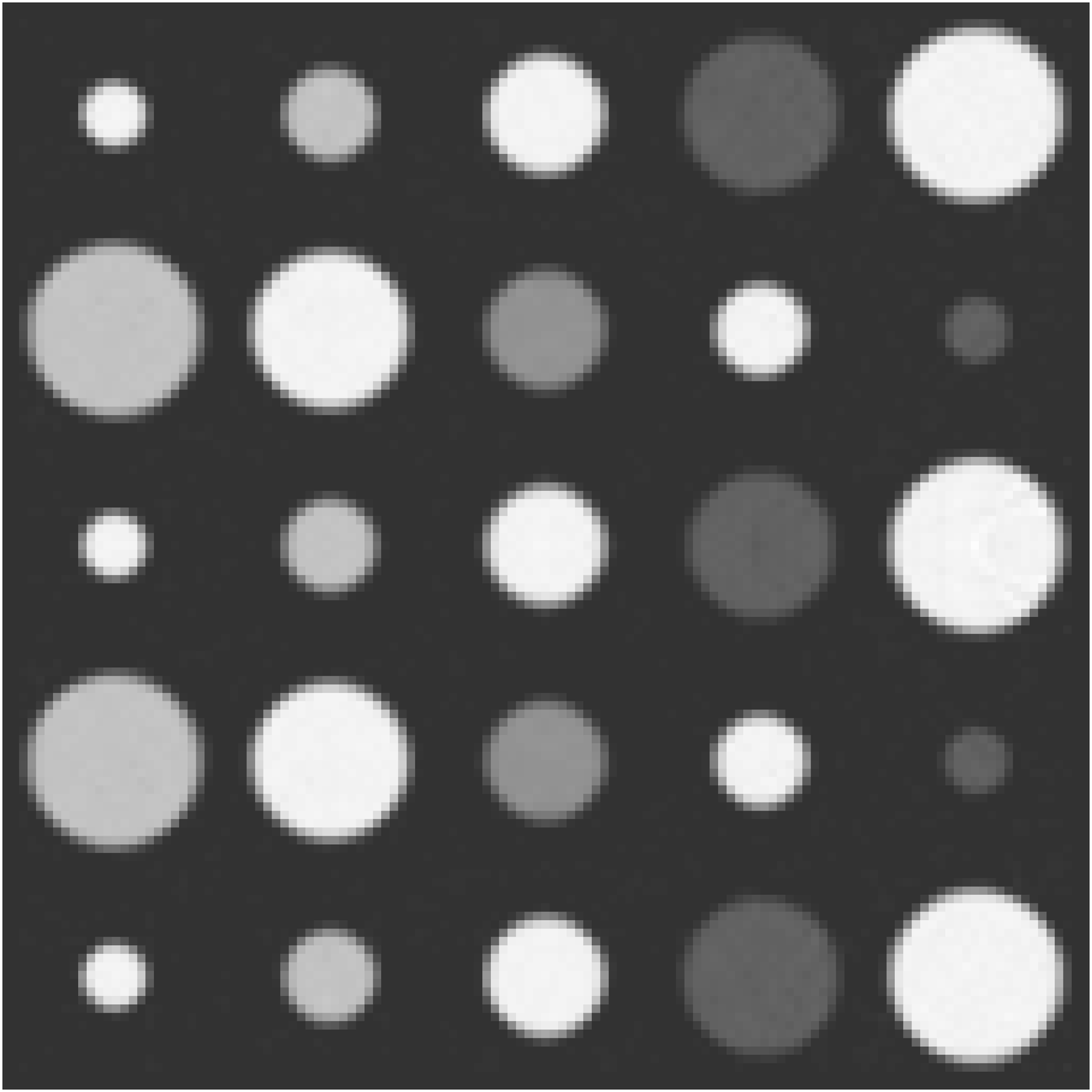} }
\end{center}
\caption{Reconstruction in 3D from spherical means centered on a surface of a
cube (a) phantom (b) reconstructed function (c) reconstruction in a presence
of a source outside of the cube}%
\label{F:3drec}%
\end{figure}

An example illustrating reconstruction of a function from its spherical means
with centers lying on the surface of a cube is presented in
Figure~\ref{F:3drec}. As a phantom we used a set of 25 characteristic
functions of balls of various radii supported within the cube $(-1,1)\times
(-1,1)\times(-1,1)$. The centers of all the balls lied in the plane $x_{3}=0.$
The cross section of our phantom by the latter plane is shown in
Figure~\ref{F:3drec}(a). The part (b) of this figure demonstrates the cross
section
of the image reconstructed on the grid of size $129\times129\times129$ using
formula~(\ref{E:3dmeaninversion}). The simulated detectors (i.e. the centers
of the integration spheres) were placed on the uniform Cartesian
$129\times129$ grids on each of the faces of the cube; there were 257
integration spheres for each detector with the radii varying uniformly from
$0$ to $2\sqrt{3}.$ Our algorithm was based on a straightforward
discretization of equation (\ref{E:3dmeaninversion}), where finite differences
were used to compute the derivatives. 

In order to illustrate the reconstruction in the presence of exterior sources
we added to the phantom shown in Figure~\ref{F:3drec}(a) a ball of radius 0.08
located at the position (1.1,0,0), and computed the spherical means
corresponding to this extended set of sources. The size and location of this
additional ball were chosen so that all the spherical means used for the
reconstruction still vanished when the radii exceeded $2\sqrt{3}.$ The result
of the reconstruction computed using formula~(\ref{E:3dmeaninversion}) is
shown in Figure~\ref{F:3drec}(c). The resulting image may seem
indistinguishable from the one in Figure~\ref{F:3drec}(b); in fact, the
difference is about 4\% in the $L^{\infty}$ norm. (In the absence of the
discretization error the images would exactly coincide).

\begin{figure}[t]
\begin{center}
\subfigure[]{\includegraphics[width=1.8in,height=1.8in]{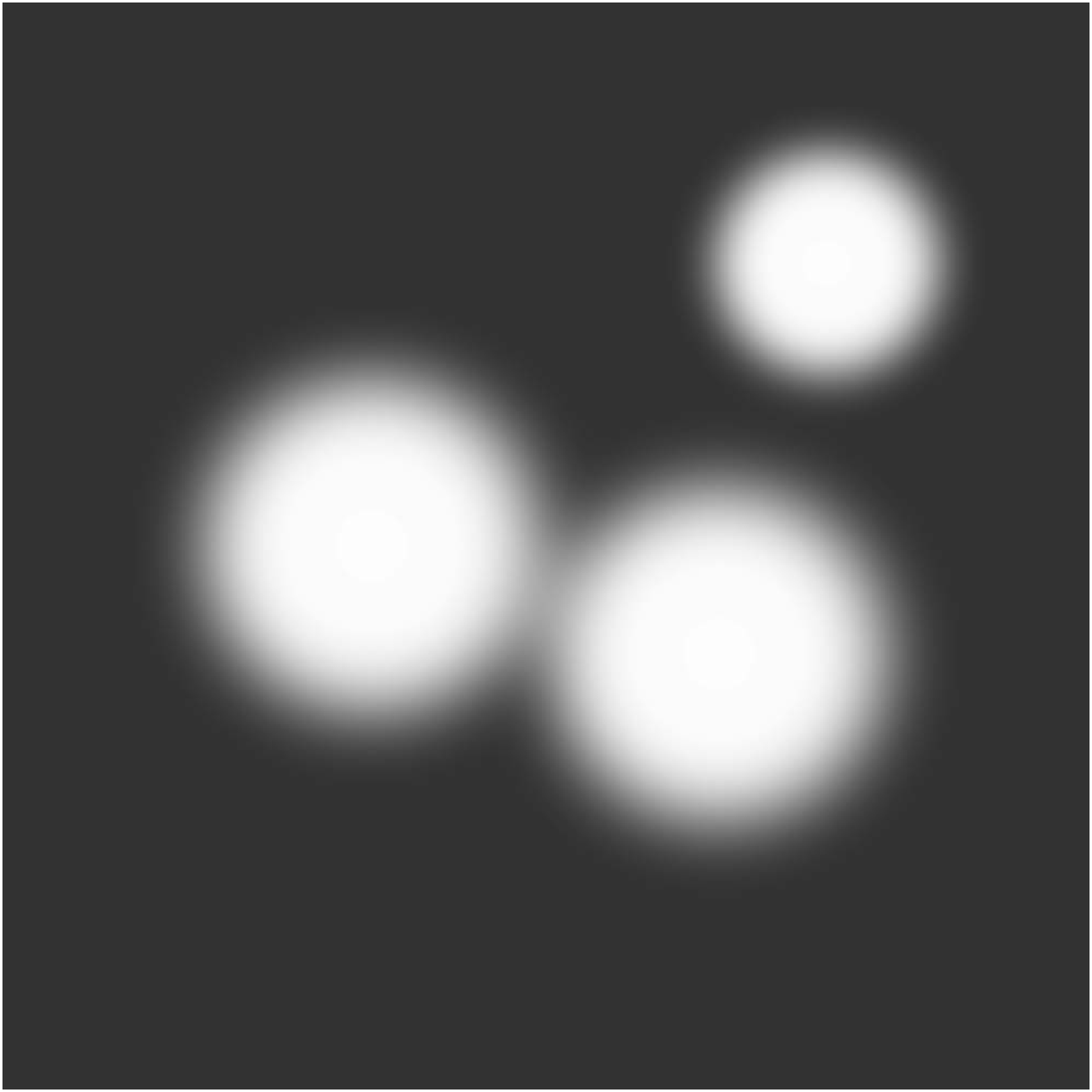} }
\subfigure[]{\includegraphics[width=1.8in,height=1.8in]{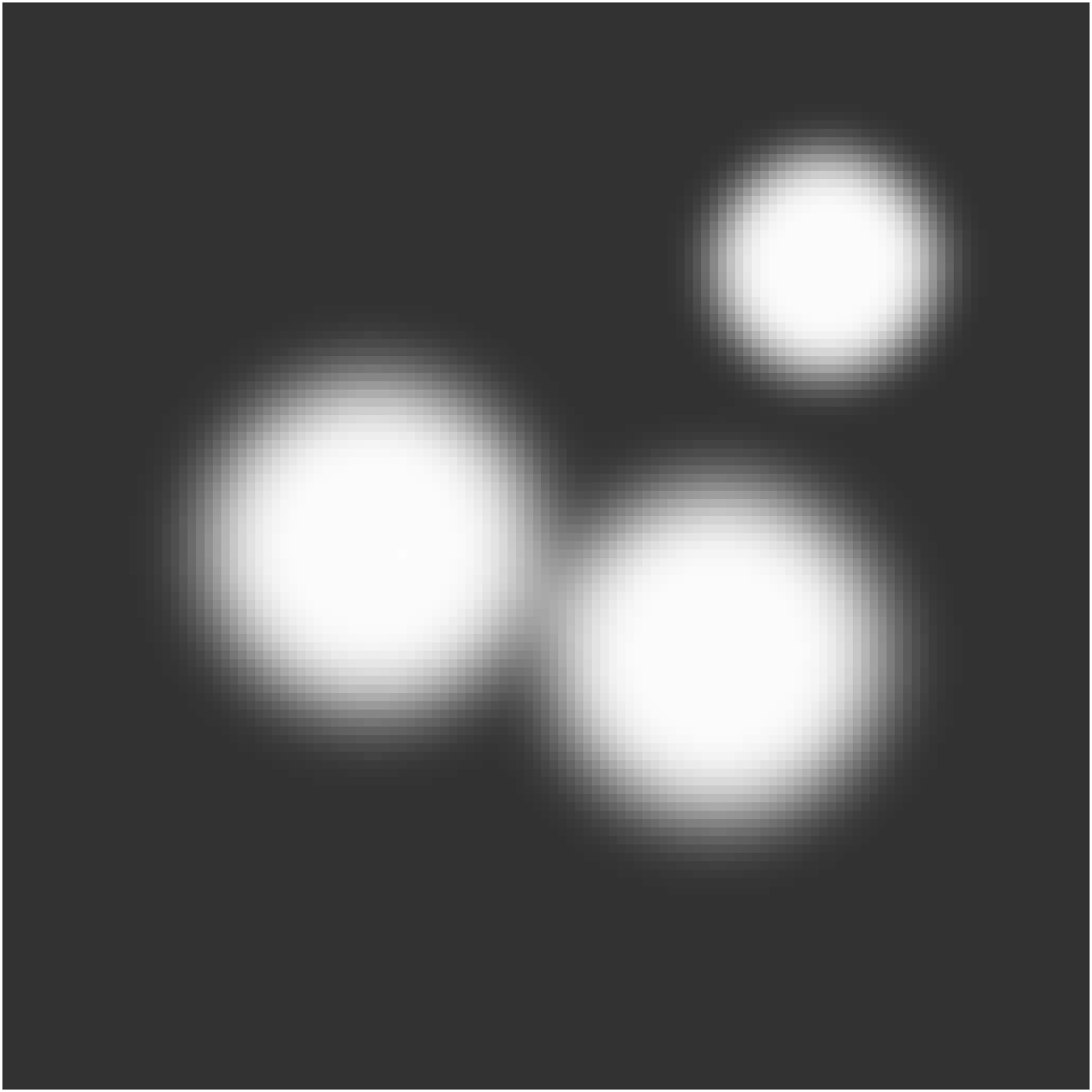} }
\end{center}
\caption{Reconstruction in 2D from circular means centered on the perimeter of
a square (a) phantom (b) reconstruction }%
\label{F:2drec}%
\end{figure}

In Figure~\ref{F:2drec} we demonstrate the reconstruction of a function in 2D
from its circular means centered on the boundary of a square $(-1,1)\times
(-1,1)$ using a modified formula~(\ref{E:2dmean1}), with the integration
restricted to the subset $\Phi^{\mathrm{trunc}}$ of $\cup_{j=1}^{4}\Phi_{j}$
contained within the disk of radius $3\sqrt{2}$ centered at the origin. (All
the points $y$ satisfying the condition $\mathrm{dist}(y,\Omega)\leq
\mathrm{diam}\Omega$ are contained in $\Phi^{\mathrm{trunc}}).$ The main goal
of this experiment was to evaluate the effect of truncating the integration
region. To this end we took care to eliminate other sources of error. In
particular, as a phantom we chose a smooth function shown in
Figure~\ref{F:2drec}(a). Moreover, our implementation of formula
~(\ref{E:2dmean1}) used a spectrally accurate algorithm for computing the
Hilbert transform and a higher order polynomial interpolation during the
backprojection. (The Hilbert transform arises since fraction $r/(r^{2}-s^{2})$
in \ref{E:2dmean1} can be re-written as $0.5/(r-s)+0.5/(r+s).)$ The
reconstructed function is shown in In Figure~\ref{F:2drec}(b). The image is
almost identical to the phantom; the relative reconstruction error in this
example equals $7.4\cdot10^{-3}$ in the $L^{\infty}$ norm. In other words, the
truncation error is negligible for most practical applications, if all the
points $y$ satisfying the condition $\mathrm{dist}(y,\Omega)\leq
\mathrm{diam}\Omega$ are included in the truncated domain.

\section{Explicit inversion formulae for some other domains\label{S:others}}

\subsection{Re-visiting the replication procedure for a cube and a square}

The inversion formulae introduced in this paper for a cube (in 3D) and a
square (in 2D) are based on the explicit representations of a solution of the
wave equation by double-layer potentials. In general, double layer potentials
are frequently used to solve Dirichlet problems for various PDE's in the
interior and/or exterior of bounded regions. Usually such formulations lead to
boundary integral equations for the potential densities that have to be solved
numerically. However, as explained in the previous sections, if the detectors
(centers of spherical means) are located on the boundary of a cube (or, in 2D,
a square), the solution of the wave equation can be represented by the double
layer potentials supported on certain periodic planes, with certain periodic
densities. Due to the symmetries in these geometries, the potentials vanish on
the planes, and the jump conditions relate the densities directly to the known
boundary values, thus providing the desired explicit reconstruction formulae.

This approach can be extended to obtain explicit reconstruction formulae for
the measurements made from the boundaries of certain polygons (in 2D) and
polyhedra in (3D) possessing sufficient symmetries. These formulae will be
similar to the present formulae for the square and the cube, but the planes
supporting the potentials and the replication procedure defining the densities
will be different in each case.

Before extending the inversion formulae to these new domains, let us provide
an alternative description of the replication procedure for a square and a
cube -- such that it would be easier to generalize to other geometries. We
start with a square, in 2D. Our goal is, again, to find the solution
$u^{T}(x,t)$ of the initial-bouindary-value problem (\ref{E:truncated2dall})
provided the boundary values of $u$ are known. To this end, let us construct a
tesselation of the plane $\mathbb{R}^{2}$ by squares of the size $2a\times2a,$
one of which coincide with $\Omega=(-a,a)\times(-a,a).$ In each square we
define $u^{\mathrm{repl}}$ to equal $\pm u^{T}$ (in some local coordinates).
The signs are chosen so as to be different in any two adjacent squares. The
local coordinates are defined in such a way that the function in any two
adjacent squares is odd with respect to the joint side of these two squares.
The resulting function $u^{\mathrm{repl}}(x,t)$ is continuous in $x$
everywhere except the boundary of each square. At each point of the boundary
it has a jump equal to a doubled negative limiting value of $u^{\mathrm{repl}%
}(x,t)$ as $x$ approaches the boundary from the inside. Now we define the
double layer potential supported on all lines containing the boundary of the
squares, with the density equal to the jump in $u^{\mathrm{repl}}.$ This new
procedure leads to the same double layer potentials as were described in
Section~\ref{S:2drev}, with the densities equal $-2P_{j}^{\mathrm{repl}}%
(x_{2},t)\eta(T,t)$ on each family of lines $\Phi_{j}$. The rest of the
reconstructing procedure is the same.

Similarly, in order to obtain the double layer potential for a solution of the
wave equation in a cube $\Omega=(-a,a)\times(-a,a)\times(-a,a)$, we tessellate
the plane with shifted versions of $\Omega$. In each shifted cube we define
$u^{\mathrm{repl}}(x,t)$ as a reflected version of $u(x,t)$ with a sign
opposite to the sign used in the neighboring cubes. The local coordinates in
each cube are defined in such a way that $u^{\mathrm{repl}}(x,t)$ is odd with
respect to each side of each cube. Now we define the double layer potentials
supported on all planes that contain all faces of all cubes, with the density
equal to the jump of $u^{\mathrm{repl}}(x,t)$ across the cube faces. This new
procedure produces the same set of the double layer potentials as was utilized
in the Section~\ref{S:3drev}, and it leads to the same reconstruction
formulae. Although the double layer potential described above is
supported on an unbounded set, only a bounded subset of it is actually
used to reconstruct the image within $\Omega$, due to the finite
speed of sound.

One obvious extension of the present inversion formulae is to rectangular
domains in 2D and to cuboids in 3D. All the symmetries needed for the double
layer potentials to form an explicit solution remain in place in these cases,
and all the reconstruction formulae (for the wave equation and for the
spherical/circular means) remain the same, with a proper re-definition of
integration lines and planes.

\subsection{Inversion formulae for certain triangular domains in 2D}

\begin{figure}[t]
\begin{center}
\subfigure[]{\includegraphics[width=1.8in,height=1.8in]{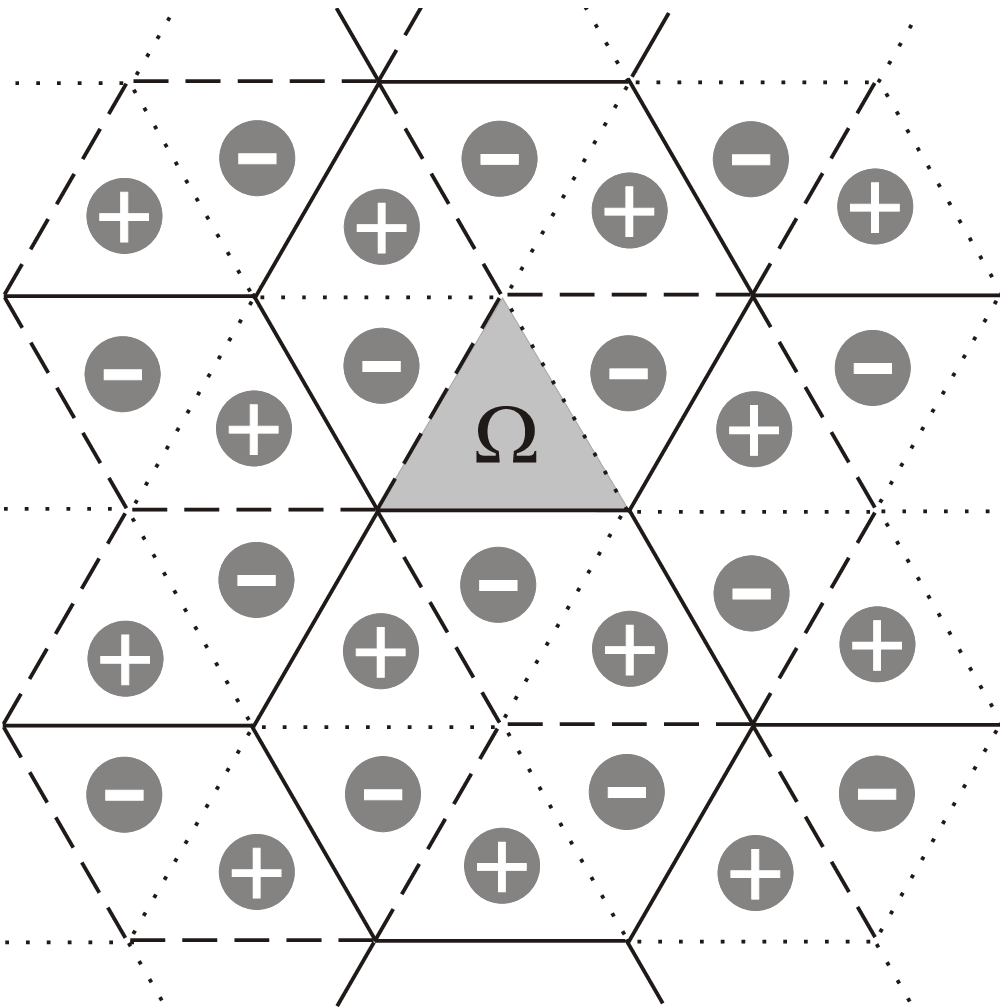} }
\subfigure[]{\includegraphics[width=1.8in,height=1.8in]{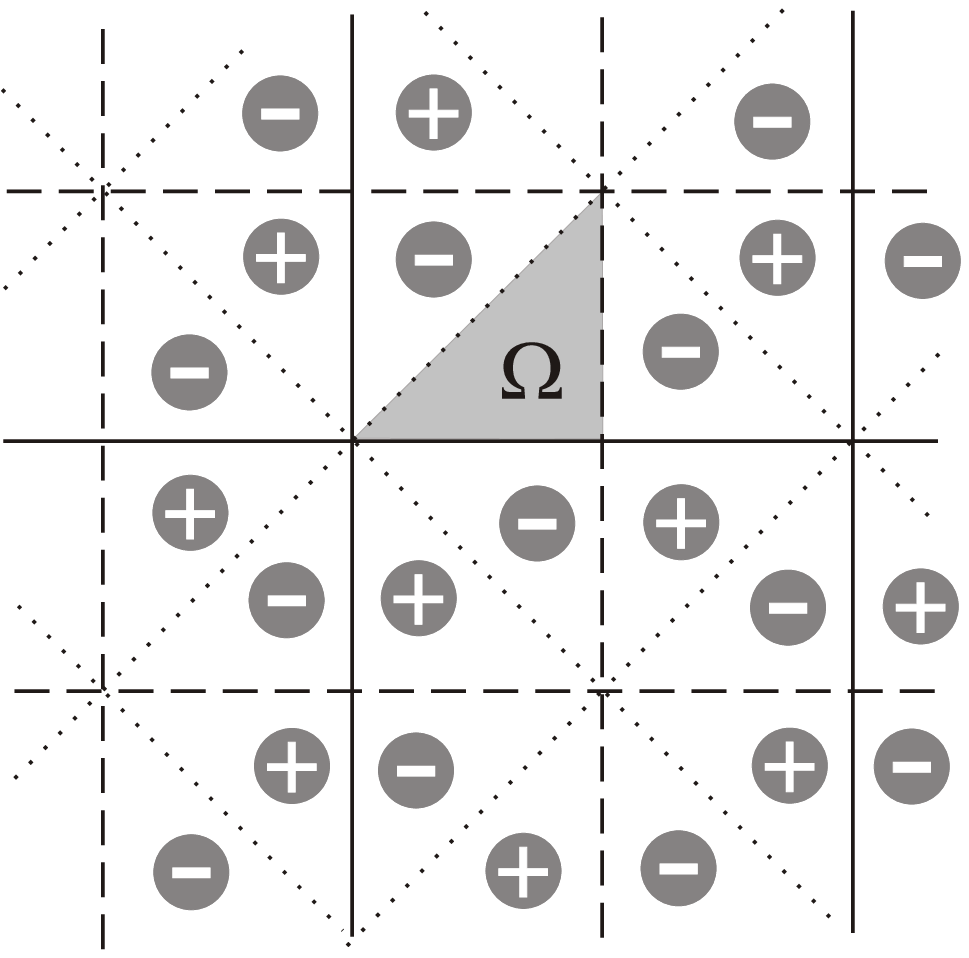}  }
\subfigure[]{\includegraphics[width=1.8in,height=1.8in]{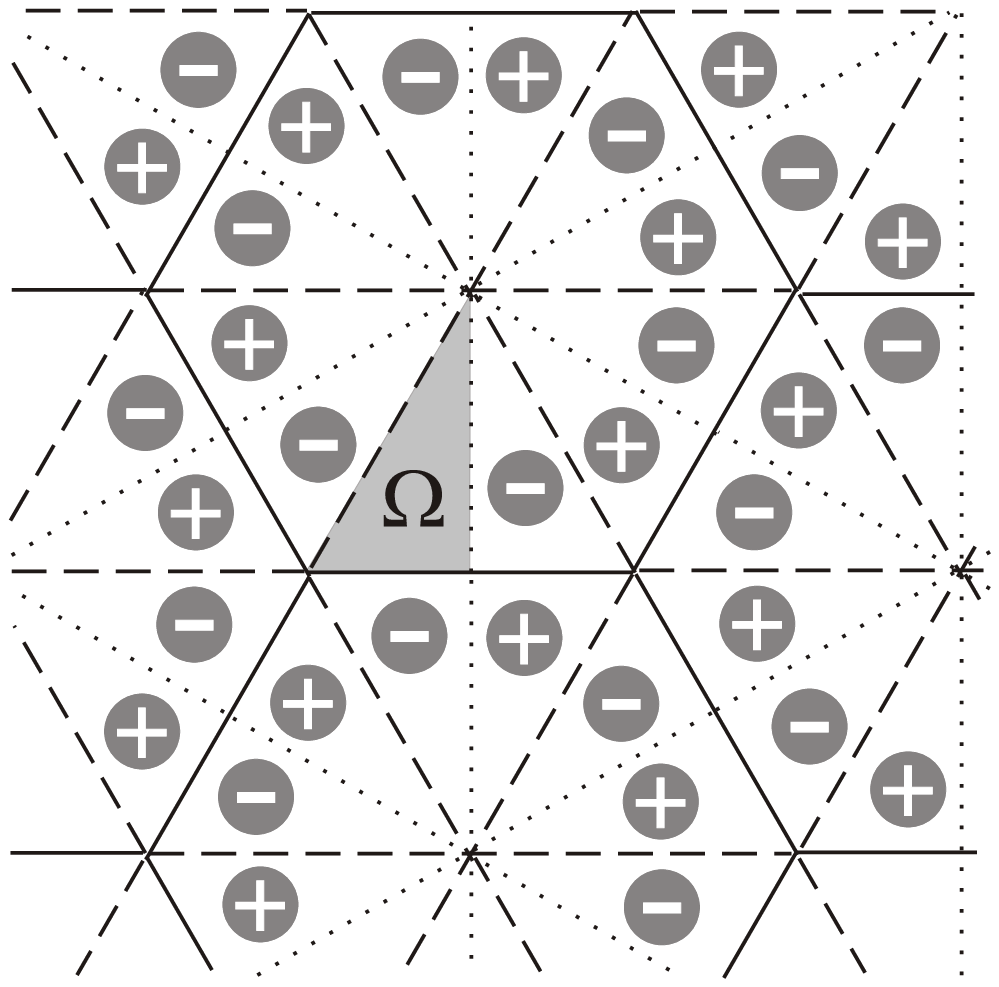}  }
\end{center}
\caption{The tesselation schemes for some triangular domains (a) equilateral
triangle (b) right isosceles triangle (c) triangle with angles $\frac{\pi}%
{2},\frac{\pi}{3}$ and $\frac{\pi}{6}$}%
\label{F:2drepl}%
\end{figure}

In 2D, the inversion formulae presented in Sections~\ref{S:2drev} and
\ref{S:2dmeans} can be generalized to certain triangular domains. Namely, a
equilateral triangle, a right isosceles triangle, and a triangle with the
angles $\frac{\pi}{2},\frac{\pi}{3}$ and $\frac{\pi}{6}$ can be used to
tessellate the plane while preserving the necessary symmetries. The
corresponding tesselations are shown in Figure~\ref{F:2drepl}.
In each of the three cases shown in the figure, the sides of the
triangular fundamental domain $\Omega$ are marked by different types of lines
(solid, dashed, and dotted ones) so that the orientation of the local
coordinates in each of the triangles is clear from the picture. Signs "+" and
"-" indicate the choice of $u$ or $-u$ in each triangle. Clearly, the
resulting $u^{\mathrm{repl}}(x,t)$ is odd with respect to all lines containing
all sides of all the triangles. Therefore, if one defines the double layer
potential with the density supported on these lines and equal to the jump of
$u^{\mathrm{repl}}(x,t)$ across the corresponding face, the value of this
potential will vanish on all the triangle sides, and, due to the jump
condition, this potential will solve the initial-boundary-value problem for
the wave equation in $\Omega.$ Therefore, the analog of the
formula~(\ref{E:finalwave2D}) is valid with the corresponding re-definition of
$P^{\mathrm{repl}}$ and with $\Omega$ being one of the above mentioned triangles.

This, in turn, allows us to obtain the analog of the inversion formulae
(\ref{E:2dmean1}) - (\ref{E:2dmean3}) for the reconstruction from circular
means, by defining $m^{\mathrm{repl}}$ in a similar fashion. As it is the case
for the square domain, the latter formulae are exact if the integration is
done over the unbounded graph. However, due to the fast decrease of the
expression in the parentheses in~(\ref{E:2dmean1}) one can hope that the
integration over the region containing all points $y$ satisfying the condition
$\mathrm{dist}(y,\Omega)\leq\mathrm{diam}\Omega$ will produce an accurate
approximation to $f(x).$ \begin{figure}[h]
\begin{center}
\includegraphics[width=6in,height=1.2in]{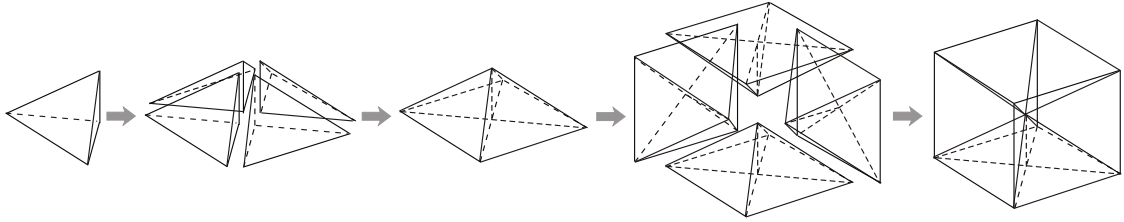}
\end{center}
\caption{Extending $u(x,t)$ by odd reflections to a cube}%
\label{F:pyramid}%
\end{figure}

\subsection{Inversion formulae for certain polyhedral domains in 3D}

In addition to cuboids, the techniques of the previous sections can also be
used to find explicit solutions and inversion formulae in 3D for some other
domains. In particular, one can derive explicit inversion formulae for
reconstructing a function from its spherical means centered on the surface of
a right prism of height $H$, whose base is one of the triangles discussed in
the previous section. As before, it may be easier to start with finding the
explicit solution to the wave equation $u(x,t)$ whose values on the boundary
of the prism are known. In order to obtain such a solution, we replicate
$u(x,t)$ in the directions parallel to the base of the prism, according to the
schemes in Figure~\ref{F:2drepl}. Then we replicate the resulting slab of
height $H$ in the odd fashion with respect to the planes containing the bases
of the prism, until $u(x,t)$ is extended to function $u^{\mathrm{repl}}(x,t)$
defined in the whole $\mathbb{R}^{3}$ in $x.$ Further, we define the double
layer potential with the density supported on planes containing the faces of
all the replicated prisms, and equal to the jumps of $u(x,t)$ across the
replicated faces. This double layer potential vanishes on all the faces, and,
due to the jump condition, provides a solution to the initial-boundary-value
problem for the wave equation in 3D. This results in an explicit formula
analogous to~(\ref{E:3dwaveinversion}). In order to obtain the analog of the
formula~(\ref{E:3dmeaninversion}) it is enough to define $M^{\mathrm{repl}}$
according to the replication procedure described above. Due to the Huygens
principle, the integration is actually performed over the subset of the
replicated faces satisfying the condition $|x-y|<\mathrm{diam}\Omega;$ both
formulae (for the spherical means and for the wave equation) are exact in this case.

Finally, the present method can be extended to the triangular pyramid $\Omega$
whose vertices have coordinates $(0,0,0),(a,0,0),(0,a,0),$ and $(0,0,a).$
(This region can also be described as a pyramid whose side faces are equal
right isosceles triangles.) Again, starting with the solution $u(x,t)$ of the
wave equation in $\Omega,$ we define $u^{\mathrm{repl}}(x,t)$ as follows.
First, we reflect $u(x,t)$ in odd fashion with respect to the vertical faces
of $\Omega$ to extend $u$ to a rectangular pyramid (see Figure~\ref{F:pyramid}%
). Then we reflect (oddly) the pyramid with respect to its triangular faces to
obtain six pyramids (only four are shown in the figure) to form a cube. Then
the cube is used to tessellate the whole space. As before, we define the
double layer potential by the jumps of the resulting function
$u^{\mathrm{repl}}(x,t)$ across all the replicated faces. The symmetries of
such an object guarantee that the potential vanishes on all the faces, and the
jump conditions yield the required limiting values as one approaches the
boundary of $\Omega$ from the inside. Thus, one obtains the exact solution to
the initial-boundary-value problem corresponding to the time reversal of the
wave equation. By implementing a similar replication procedure to define
$M^{\mathrm{repl}},$ one obtains an analog of the inversion
formula~(\ref{E:3dmeaninversion}) for reconstructing a function from its
spherical means centered on the surface of $\Omega.$

It is worth noting that, similarly to the formulae for a cube and a square,
all the formulae discussed in the present and previous sections also have the
property of being insensitive to sources lying outside of the domain
surrounded by the acquisition surface.

\section*{Conclusions}

By utilizing the double layer wave potentials we found explicit solutions to
the initial-boundary-value problems for the wave equation in certain 2D and 3D
domains. In problems of TAT/PAT this yields explicit expressions for the
result of the time reversal, and thus explicitly reconstructs the initial
pressure distribution $f(x).$ Further, by formulating the problem in terms of
spherical (circular) means, we obtain from the solutions to the wave equations
explicit closed-form inversion formulae of filtration/backprojection type for
the corresponding domains. In 2D, we presented inversion formulae for a square
(or a rectangle), and we outlined their generalization for an equilateral
triangle, right isosceles triangle, and a triangle with the angles $\frac{\pi
}{2},\frac{\pi}{3}$ and $\frac{\pi}{6}$. In 3D, we derived formulae for a cube
(with immediate extension to a cuboid), and we discussed their generalization
to right prisms whose bases are one of the three above-mentioned triangles,
and to a pyramid whose side faces are equal right isosceles triangles.

In all cases the formulae for the inversion of the wave equation data are in
the form of the double layer potentials supported on a series of planes, with
the densities whose values equal to the known boundary values replicated
according to the procedures described in the text. The formulae for inverting
the spherical (or circular) means easily follow from the ones for the wave
data. In the 3D case, due to the Huygens principle the integration in all the
formulae is restricted to a bounded subset of the planes. In 2D, in order to
obtain the exact reconstruction one has to integrate over an unbounded domain.
However, our numerical experiments suggest that for a sufficiently accurate
reconstruction it is enough to integrate over a bounded subset containing all
points
$y$ satisfying the condition $\mathrm{dist}(y,\Omega)\leq\mathrm{diam}\Omega.$

Finally, all the presented inversion formulae remain exact within the domain
even in the presence of sources located outside the domain. (In order to
observe this phenomenon in 3D one may need to extend the measurement time (or
the radii of the spherical means) to make sure that the signal vanishes after
the observation time). Interestingly, although the time reversal and some
series methods also have this property, the present formulae are the only
known explicit FBP-type inversion formulae that exhibit such behavior.

\paragraph*{Acknowledgements}

The author would like to thank Prof. P. Kuchment for helpful discussions, and
to gratefully acknowledge support by the NSF through the grant DMS-090824.


\begin{thebibliography}{99}
\bibitem {FR2}D.~Finch and Rakesh, Recovering a function from its spherical
mean values in two and three dimensions. In \cite{WangCRC}, 77--88.
\end{thebibliography}
\end{document}